\documentclass[11pt,oneside]{amsart}
\usepackage{amsmath,amssymb,graphicx,subfigure,psfrag}

\newtheorem{thm}{Theorem}%[section]
\newtheorem{lem}[thm]{Lemma}
\newtheorem{cor}[thm]{Corollary}

\theoremstyle{definition}
\newtheorem{defn}[thm]{Definition}

\newtheorem{conj}[thm]{Conjecture}

\newcommand{\bi}{\begin{itemize}}
\newcommand{\ei}{\end{itemize}}
\newcommand{\be}{\begin{enumerate}}
\newcommand{\ee}{\end{enumerate}}
\newcommand{\bc}{\begin{center}}
\newcommand{\ec}{\end{center}}
\newcommand{\bt}{\begin{tabular}}
\newcommand{\et}{\end{tabular}}

\newcommand{\Z}{\mathbb Z}
\newcommand{\N}{\mathbb N}

\newcommand{\M}{\mathcal M}\newcommand{\G}{\mathcal G}

\begin{document}
%%%%%%%%%%%%%%%%

\title{$G$-automata, counter languages and the
 Chomsky hierarchy}

\author{Murray Elder}
 \thanks{Supported by EPSRC grant GR/S53503/01}

\address{School~of
Mathematics and Statistics, University of St~Andrews, KY16~9SS Scotland}
\email{murrayelder@gmail.com, http://me.id.au/}

 \date{\today}
\begin{abstract}
We consider how the languages of $G$-automata compare with other formal language classes. We prove that if the word problem of $G$ is accepted by a machine in the class $\M$ then the language of any $G$-automaton is in the class $\M$. It follows that the so called {\em counter languages} (languages of  $\Z^n$-automata) are context-sensitive, and further that counter languages are indexed if and only if  the word problem for $\Z^n$ is indexed.
\end{abstract}

  

\maketitle

\section{Introduction}
In this article we 
 compare the languages of $G$-automata, which include the set of {\em counter languages}, with the formal language classes of context-sensitive, indexed,  context-free and regular. We prove in Theorem \ref{thm:Main} that if the word problem of $G$ is accepted by a machine in the class $\M$ then the language of any $G$-automaton is in the class $\M$. It follows that the  counter languages  are context-sensitive. Moreover it follows that counter languages are indexed if and only if  the word problem for $\Z^n$ is indexed. 
 
  The article is organized as follows. In Section \ref{sec:defns} we define $G$-automata, 
  linear-bounded automata, nested-stack, stack, and pushdown automata, and the word problem for a
   finitely generated group. In  Section \ref{sec:Thm} we prove the main theorem,
and give the corollary that counter languages are indexed if and only 
if the word problem for $\Z^n$ is indexed for all $n$.

\section{Definitions}\label{sec:defns}
If $G$ is a group with generating set $\G$, we say two words $u,v$
are equal in the group, or $u=_Gv$, if they represent the same group
element. We say $u$ and $v$ are identical if the are equal in the free
monoid, that is, they are equal in $\G^*$.

\begin{defn}[Word problem]
The {\em word problem} of a group $G$ with respect to a finite generating set $\G$
is the set of words $\{w\; |\; w=_G1\}$. Note that the word problem for a finitely generated group is a language over a finite alphabet.
\end{defn}

\begin{defn}[$G$-automaton]\label{def:G-aut} 
Let $G$ be a group and $\Sigma$ a finite set.  A (non-deterministic)
 {\em $G$-automaton} $A_G$ over $\Sigma$ is a finite directed graph
 with a distinguished {\em start vertex} $q_0$, some distinguished
 {\em accept vertices}, and with edges labeled by $(\Sigma^{\pm
 1}\cup\{\epsilon\})\times G$.  If $p$ is a path in $A_G$, the element
 of $(\Sigma^{\pm 1})$ which is the first component of the label of
 $p$ is denoted by $w(p)$, and the element of $G$ which is the second
 component of the label of $p$ is denoted $g(p)$. If $p$ is the empty
 path, $g(p)$ is the identity element of $G$ and $w(p)$ is the empty
 word.  $A_G$ is said to {\em accept} a word $w\in (\Sigma^{\pm 1})$
 if there is a path $p$ from the start vertex to some accept vertex
 such that $w(p)=w$ and $g(p)=_G 1$.
 \end{defn}

\begin{defn}[Counter language]
A language is {\em $k$-counter} if it is accepted by some
$\Z^k$-automaton. We call the (standard) generators of $\Z^k$ 
{\em counters}. A language is {\em counter} if it is $k$-counter for
some $k\geq 1$.
\end{defn}

These definitions are due to 
Mitrana and Striebe \cite{Mitrana}. Note that in these counter automata, the values of the counters is not accessible until the final accept/fail state. For this reason they are sometimes called {\em blind}.
Elston and Ostheimer \cite{EO} proved that the word problem of $G$ is deterministic counter with an extra "inverse" property  if and only if $G$ is virtually abelian. Recently  Kambites \cite{Kambites} has shown that the inverse property restriction can be removed from this theorem.

It is easy to see that the word problem for $G$ is accepted by a $G$-automaton.
\begin{lem}\label{lem:WPG}
The word problem for a finitely generated group $G$ is accepted by a deterministic $G$-automaton.
\end{lem}
\textit{Proof}:
Construct a $G$-automaton with one state and a directed loop labeled by $(g,g)$ for each generator $g$. The state is both start and accept. A word in the generators is accepted by this automaton if and only if it represents the identity, by definition.
\hfill$\Box$

Recall the definitions of the formal language classes of recursively enumerable, decidable, context-sensitive, indexed, stack, context-free, and regular. Each of these can be defined as the languages of some type of restricted  Turing machine as follows.

Consider a machine consisting of finite {\em alphabet} $\Sigma$, a finite {\em tape alphabet} $\Gamma$, a {\em finite state control} and an infinite {\em work tape}, which operates as follows. The finite state control is a finite graph with a specified {\em start node}, some specified {\em accept nodes}, and edges labeled by an alphabet letter and an instruction for the work tape. The instructions in general are are the form \texttt{read, write, move left, move right} and one reads/writes letters from $\Gamma$ on the tape. The tape starts out blank.

One inputs a finite string in $\Sigma^*$ one letter at a time, read from left to right. For each letter $x\in \Sigma$, the finite state control performs  some instructions on the work tape corresponding to an edge who label is $(x, \mathtt{instructions})$, and moves to the target node of the edge. One starts at the start node, and {\em accepts} the string if there is some path from the start node to an accept node labeled by the letters of the string. The {\em language} of a machine is the set of strings in $\Sigma^*$ which the machine accepts. If the finite state control includes edges of the form $ (\epsilon, \mathtt{instructions})$
then one can work on the work tape without reading input, and if a machine has such edges, or two edges with the same letter in the first coordinate of the edge label, then such machines (and their languages) are called {\em non-deterministic}.

If we allow no further restrictions on how this machine operates, then we have a {\em Turing machine}.  By placing increasingly strict restrictions on the machine, we obtain a hierarchy of languages corresponding to the machines.
If the Turing machine halts on accepted strings then the language is {\em recursively enumerable} or {\em r.e.}, and if it halts both accepted and rejected strings the language is {\em decidable}.
If we restrict the number of squares of tape that can be used to be a constant multiple of the length of the input string, then we obtain a {\em linear bounded automaton}, and the languages of these are called {\em context-sensitive}.
If we make the tape act as a {\em nested stack} (see \cite{GilmShap}) then the language of such a machine is called {\em indexed}. 

If the tape is a stack (first in last out) where the pointer may read but not write on any square, then  the machine and its languages are  called {\em stack}.
If the tape is a stack such that the pointer can only read the top square, we have a {\em pushdown automaton}, the languages of which are {\em context-free}.
Finally, if we remove the tape altogether we are left with just the finite state control, which we call a {\em finite state automaton}, languages of which are {\em regular}.

%The language of a Turing machine which has a constant bounded number of moves on the tape for each input letter, then the machine and its languages are  called {\em real-time}. Real-time languages do not strictly fit into the above hierarchy. (Example-Derek)

For more precise definitions
see \cite{GilmShap} for nested stack automata, \cite{HoltRover} for stack, \cite{Sipser} for regular and context-free  and \cite{HU} for these plus linear bounded automata.
Two good survey articles are \cite{Hairdressing} and \cite{Gilman}.

\section{Main theorem}\label{sec:Thm}

\begin{thm}\label{thm:Main}
Let $\M$ be a formal language class:
(regular, context-free, stack, indexed, context-sensitive, decidable, r.e.) and let $G$ be a finitely generated group. The word problem for $G$ is in $\M$ if and only if the language of
every $G$-automaton is in $\M$.
\end{thm}
\textit{Proof}:
By Lemma \ref{lem:WPG} the word problem of $G$ is accepted by a $G$-automaton, so one direction is done.

Let $L$ be a language over an alphabet $\Sigma$ accepted by a $G$-automaton $P$. Fix a finite generating set $\G$ for $G$ which includes all elements of $G$ that are the second coordinate of an edge label in $P$. Let $N$ be an $\M$-automaton which accepts the word problem for $G$ with respect to this generating set.

Construct a machine $M$ from the class $\M$ to accept $L$ as follows. 
The states of $M$ are of the form $(p_i,q_j)$ where $p_i$ is a state of $P$ and $q_i$ is a state of $N$. The start state is $(p_S,q_S)$, and accept states are $(p_A,q_A)$ where $p_S,p_A,q_S,q_A$ are
start and accept states in $P$ and $N$.

The transitions are defined as follows. For each edge $(x,g)$ in $P$ from $p$ to $p'$, and for each edge 
$(g,\mathtt{instruction})$ in $N$ from $q$ to $q'$,  add an edge from $(p,q)$ to $(p',q')$ in $M$ labeled 
$(x,\mathtt{instruction})$.

Then $M$ accepts a string in $\Sigma^*$ if there is a path in $M$ corresponding to paths in $P$ from the start state to an accept state such that the  labels of the second coordinate of the edges give a word $w$ in $\G^*$, and a path in $N$ from a start state to an accept state labeled by $w$ which evaluates to the identity of $G$ and respects the tape instructions.

Note that $M$ is deterministic if both $N$ and $P$ are deterministic and no state in $P$ has two outgoing edges with the same group element in the first coordinate of the edge labels. 
\hfill{$\Box$}

It is easy to see that if $G$ is a finite group, then the construction describes a finite state automaton (we don't need any tape).  
Herbst \cite{Herbst} showed that the word problem for a group is regular if and only if the group is finite.
If $G$ is virtually free then Muller and Schupp \cite{MS} proved that $G$ has a context-free word problem, so the language of 
every $G$-automaton for $G$ virtually free is context-free.

Since the word problem for $G$ virtually abelian is context-sensitive, we get the following corollary.
\begin{cor}
Counter languages are context-sensitive.
\end{cor}
 More generally,
Gersten, Holt and Riley  \cite[Corollary B.2]{GHR} have shown that  finitely generated nilpotent groups of class $c$  have context-sensitive word problems, so the 
language of every $G$-automaton for $G$ finitely generated  nilpotent is context-sensitive.

\begin{figure}[ht!]
  \bc
     \includegraphics[width=12.5cm]{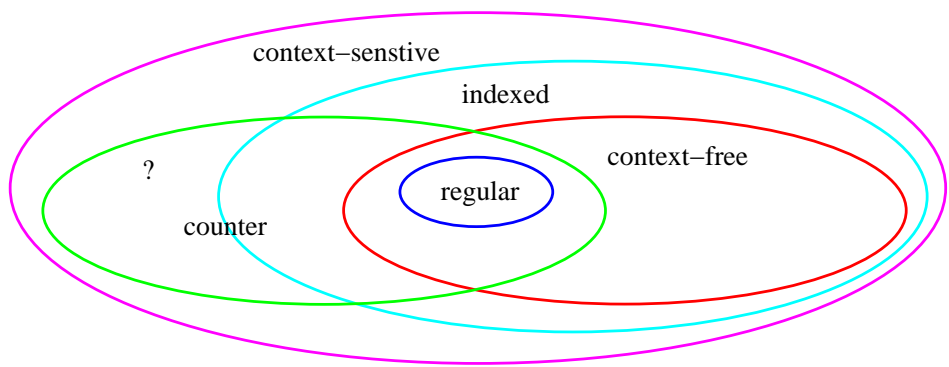}
  \ec
  \caption{How counter languages fit into the hierarchy}
  \label{fig:sets1}
\end{figure}

 Figure \ref{fig:sets1} shows how counter languages fit into the hierarchy.

\begin{cor}
The word problem for $\Z^n$ is indexed for all $n\geq 1$ if and only if every counter language is indexed.
\end{cor}

\textit{Proof:}
By Theorem \ref{thm:Main} if the word problem for $\Z^n$ is indexed then the language of every $\Z^n$-automaton is indexed, proving one direction.

Since the word problem of $\Z^n$ is an $n$-counter language then if it is not indexed then not every counter language
 is  indexed.
\hfill{$\Box$}

By Gilman and Shapiro, if the word problem of $G$ is accepted by a nested stack automaton with some extra restrictions, 
then it is virtually free. We still do not know whether word problem of $\Z^n$ is accepted by a nested stack automaton
 without extra restrictions.

\begin{conj}
The word problem for $\Z^2$ is not indexed.
\end{conj}

To motivate this conjecture and suggest a possible proof strategy, define $L_n(y,z)$ to be the set of  words in 
$\{y,z\}^{\ast}$ with $n$  $y$s and $n$ $z$s, and  consider the language 
$L=\{x^n w_n \; | \; n \in \N, w_n\in L_n(y,z) \}$. This language is the word problem $\Z^2$ with respect to the 
generators $a,b,a^{-1},b^{-1}$,  intersected with the regular language $\{(ab)^nw_n\: | \; w_n\in L_n(a^{-1},b^{-1})\}$. 
So if one could prove $L$ is not indexed then one would prove the conjecture.

\section{Acknowledgments}

Thanks to Ray Cho, Andrew Fish, Bob Gilman, Claas R{o}ver, Sarah Rees and Gretchen Ostheimer for many
 fruitful discussions about counter, context-sensitive and indexed languages and word problems.

\end{document}